\documentclass{amsart}
\usepackage{mathptmx}
\usepackage{amsmath,amsthm}
\usepackage{amscd,tikz}
\usepackage{amsfonts,enumerate}
\usepackage{amssymb}
\usepackage{graphicx,cite}
\usepackage{mathrsfs}
\usepackage[all,cmtip,2cell]{xy}

\UseTwocells
\xyoption{poly}

\newtheorem{thm}{Theorem}[section]

\theoremstyle{definition}
\newtheorem{dfn}{Definition}[section]
\theoremstyle{remark}
\newtheorem{rmk}{Remark}[section]

\title[The mathematical work of K.S.S. Nambooripad]{The mathematical work of K.S.S. Nambooripad}
\author{John Meakin}
\address{Department of Mathematics,	University of Nebraska-Lincoln\\ Lincoln, NE 68588-0323, United States of America.}
\email{jmeakin@unl.edu}
\author{P.A. Azeef Muhammed}
\address{Department of Mathematics and Natural Sciences, Prince Mohammad Bin Fahd University\\ Al Khobar 31952, Kingdom of Saudi Arabia.}
\email{aparayilajmal@pmu.edu.sa, azeefp@gmail.com}
\author{A.R. Rajan}
\address{Department of Mathematics,	University of Kerala\\ Thiruvananthapuram 695581, India.}
\email{arrunivker@yahoo.com}
\subjclass[2010]{}

\begin{document}
\maketitle
\begin{abstract}
		
We provide an overview of the mathematical work of K.S.S. Nambooripad, with a focus on his contributions to the theory of regular semigroups. In particular, we outline Nambooripad's seminal contributions to the structure theory of regular semigroups via his theory of {\em inductive groupoids}, and also via his theory of {\em cross connections}. We also provide information about outgrowths of his work in the algebraic theory of semigroups and its connections with several other fields of mathematics, in particular with the theory of operator algebras.
		
\end{abstract}
	
\section{Introduction}
	
K.S.S. Nambooripad passed away in his home in Trivandrum, India on January 4, 2020. He was an outstanding mathematician and the founder of an important school of semigroup theory in the Indian state of Kerala. Some  information about his work may be found in the tribute \cite{MeRa} that the first and third authors of this article wrote on the occasion of Nambooripad's 80th birthday in 2015. Some interesting additional information about Nambooripad's life and his impact on the development of the \text{\TeX}  software package in India may be found in Radhakrishnan's blog \cite{Rad}.  In the present article, we summarize some of his mathematical work and its impact on the development of semigroup theory.
	
Nambooripad's most influential work was primarily concerned with the study of (von Neumann) regular semigroups and their connections with other fields of mathematics, and this article will focus on some of his contributions in this area. He made seminal contributions to the study of regular semigroups through his two deep works \cite{KSS2} and \cite{KSS3}. The first of these develops a structure theory of regular semigroups via his notions of {\em biordered sets} and {\em inductive groupoids}, while the second studies regular semigroups via his theory of {\em normal categories} and {\em cross connections}. Both works have had a major impact on our understanding of regular semigroups and on the development of much subsequent and ongoing work in the field.
	
We will assume that the reader is familiar with the basic ideas and notation of semigroup theory, as presented in the books of Clifford and Preston \cite{CP} or Howie \cite{How}. We will also make use of  the basic  notions of category theory as presented in the books by MacLane \cite{mac} or Higgins \cite{Hig}. If $\mathcal{C}$ is a category, then we will use $v\mathcal{C}$ to denote the set of vertices (objects) of the category $\mathcal C$ and we will denote the set of morphisms in $\mathcal C$ from $v$ to $w$ by $\mathcal{C}(v,w)$. For most categories we will denote the composition of a morphism $f \in \mathcal{C}(v,w)$ and $g \in \mathcal{C}(w,z)$ by $fg$, but this will be denoted by $gf$ in the section on Fredholm operators in Section 6 of the paper.
	
Section 2 of this paper is based on Nambooripad's memoir \cite{KSS2}. It outlines his concept of biordered sets and inductive groupoids and their role in the structure theory of regular semigroups. In Section 3, we discuss some of the outgrowths of Nambooripad's inductive groupoid approach to the structure of regular semigroups and its impact on subsequent and ongoing work in the field. In particular, we discuss his concept of the natural partial order on a  regular semigroup and also his construction of regular idempotent-generated semigroups with a  given biordered set of idempotents. We also provide some references to subsequent work extending the inductive groupoid approach to classes of non-regular semigroups. Section 4 outlines Nambooripad's work on the construction of regular semigroups via his theory of normal categories and cross connections, based primarily on his manuscript \cite{KSS3}. In Section 5, we discuss the connection between Nambooripad's two manuscripts \cite{KSS2} and \cite{KSS3}. We also discuss the cross connection theory in several special cases of regular semigroups and its extension to some classes of non-regular semigroups. Section 6 outlines some of Nambooripad's work on connections between regular semigroups and other areas of mathematics, with particular emphasis on his work in operator algebras. We close the paper in Section 7 with some brief information about Nambooripad's PhD students.


\section{Regular semigroups, biordered sets and inductive groupoids}\label{secind}

In the present section, we describe how regular semigroups give rise to biordered sets and inductive groupoids and conversely how inductive groupoids give rise to regular semigroups. We also describe Nambooripad's construction of fundamental regular semigroups. Our treatment is somewhat informal, focussed on the key ideas of Nambooripad's construction: the reader is referred to Nambooripad's memoir \cite{KSS2} for full details of this aspect of his work.

\medskip

{\bf Biordered sets}

\medskip

Let $S$ be a  semigroup with set $E(S)$ of idempotents. Define two quasi-orders (that is, reflexive and transitive relations) $\mathrel{\omega} ^l$ and $\mathrel{\omega}^r$    on $E(S)$  as follows.
\[e \mathrel{\omega}^r f\text{ if }fe=e\text{ and }e\mathrel{\omega}^l f\text{ if }ef=e.\] Then $\mathscr{R}=\omega^r \cap (\omega^r)^{-1}$ and  $\mathscr{L}=\omega^l \cap (\omega^l)^{-1}$ are equivalence relations on $E(S)$ and $\omega=\omega^r \cap \omega^l$ is the natural partial order on $E(S)$.

It is easy to see that if $e$ and $f$ are
idempotents of  $S$ such that $e\mathrel{\omega}^l f$ or $f\mathrel{\omega}^l e$ or $e\mathrel{\omega}^r f$ or $f\mathrel{\omega}^r e$ then the product $ef$ is another idempotent of $S$.
Nambooripad refers to such products as {\it basic products}. Thus
the set $E(S)$ of idempotents of $S$ becomes a partial binary
algebra with respect to the basic products. Nambooripad provided an axiomatic characterization of the partial binary algebra of idempotents of a regular semigroup with respect to the basic products and referred to such a partial algebra as a  {\em regular biordered set}. We emphasize that there may be other (non-basic) products in $E(S)$ that are also idempotents. Clifford \cite{cliff1} characterized the partial algebra of idempotents of a regular semigroup $S$ with respect to the partial binary operation $e * f = ef$ defined whenever the product $ef$ is an idempotent of $S$.

Nambooripad's axiomatic characterization of $E(S)$ with respect to the basic products is intrinsic to $E(S)$ in
much the same way as the set of idempotents of an inverse semigroup
may be characterized as a (lower) semilattice.
In fact his notion of a biordered set may be viewed as a
vast generalization of the notion of a  semilattice. The role of the
meet (product) of idempotents in an inverse semigroup is served by
what Nambooripad refers to as the {\it sandwich set} $\mathcal{S}(e,f)$ of two
idempotents $e,f \in E(S)$. If $S$ is a  semigroup and $e,f \in
E(S)$, then one may define

$$\mathcal{S}(e,f) = \{h \in E(S) : he = h = fh, \, ef = ehf\}.$$

Note that if $e$ and $f$ are idempotents of a {\it regular}
semigroup $S$, then $\mathcal{S}(e,f) \neq \emptyset$. In fact if $(ef)'$ is
any inverse of $ef$ in $S$ then the element $h = f(ef)'e$ is in
$\mathcal{S} (e,f)$. (Readers familiar with ``Lallement's Lemma" will see this
as related to an observation that Lallement made in his proof that
idempotents lift under morphisms between regular semigroups
\cite{Lall}.) Nambooripad characterized the sandwich set $\mathcal{S}(e,f)$
completely in terms of basic products, namely
$$h\in \mathcal{S}(e,f) \iff h\mathrel{\omega}^l e,\: h\mathrel{\omega}^r f\text{ and if } g\mathrel{\omega}^l e,\: g\mathrel{\omega}^r f,\text{ then } eg\mathrel{\omega}^r eh\text{ and }gf \mathrel{\omega}^l hf.$$

Nambooripad defined a {\it biordered set} to be a partial binary algebra $E$ satisfying the following axioms:
for $e,f,g \in E$,
\begin{enumerate}
	\item [(B1)] $ef$ is defined if and only if $e\mathrel{\omega}^l f$ or $f\mathrel{\omega}^l e$ or $e\mathrel{\omega}^r f$ or $f\mathrel{\omega}^r e$;
	\item [(B2)] if $f\mathrel{\omega}^l e$, then $f \mathrel{\mathscr{L}} ef \mathrel{\omega} e$; if $f\mathrel{\omega}^r e$, then $f \mathrel{\mathscr{R}} fe \mathrel{\omega} e$;
	\item [(B3)] if $g\mathrel{\omega}^r f$ and $f,g\mathrel{\omega}^l e$, then $eg \mathrel{\omega}^r ef$, $gf\mathrel{\omega}^l e$ and $e(gf) = (eg)(ef)$;\\
	if $g\mathrel{\omega}^l f$ and $f,g\mathrel{\omega}^r e$, then $ge \mathrel{\omega}^l fe$, $fg\mathrel{\omega}^r e$ and $(fg)e = (fe)(ge)$;
	\item [(B4)] if $g\mathrel{\omega}^l f \mathrel{\omega}^l e$, then $fg= f(eg)$; if $g\mathrel{\omega}^r f \mathrel{\omega}^r e$, then $gf= (ge)f$;
	\item [(B5)] if $f,g\mathrel{\omega}^l e$, then $\mathcal{S}(ef,eg)= e\mathcal{S}(f,g)$; if $f,g\mathrel{\omega}^r e$, then $\mathcal{S}(fe,ge)= \mathcal{S}(f,g)e$.
\end{enumerate}

He defined a biordered set $E$ to be a {\it regular biordered set} if in addition

\begin{enumerate}
\item [(R)] $\mathcal{S}(e,f) \neq \emptyset$ for all $e,f \in E$.

\end{enumerate}

Nambooripad went on to define the notion of a {\it bimorphism} between biordered sets (and a {\it regular bimorphism} between regular biordered sets) in a natural way, thus enabling
the class of biordered sets and the class of regular biordered sets to be viewed as  categories.
The following theorem is proved in Nambooripad's memoir \cite[Theorem 1.1 and Corollary 4.15]{KSS2}.

\begin{thm}
\label{regbio}

If $S$ is a regular semigroup, then $E(S)$ is a regular biordered set relative to the basic products. The assignment $S \rightarrow E(S), \, \phi \rightarrow  {\phi}|_{E(S)}$ is a functor from the category of (regular) semigroups to the category of (regular)  biordered sets. Conversely, every regular biordered set is (biorder) isomorphic to the biordered set of some regular semigroup.

\end{thm}

While the proof of the ``direct" part of Theorem \ref{regbio} is relatively straightforward, the proof of the ``converse " part is considerably more difficult, involving Nambooripad's theory of {\em inductive groupoids}, outlined below.

Nambooripad \cite{KSS1} originally provided an equivalent axiomatic characterization of the
biordered set of a regular semigroup
$S$ as the set $E(S)$ equipped with the two quasiorders ${\omega}^r$
and ${\omega}^l$  together with translations $\tau^{r}(e)$ and $\tau^{l}(e)$
associated with these quasiorders that enable us to define the basic
products. Here

\[f\tau^{r}(e) = fe \text{ if }f \in \omega^{r}(e) \text{ and } f\tau^{l}(e) = ef \text{ if } f \in \omega^{l}(e).\]

If one drops the requirement that $\mathcal{S}(e,f) \neq \emptyset$ for each
$e,f $, one obtains the axioms for a (not necessarily regular)
biordered set. Easdown \cite{Eas1}  subsequently showed that
Nambooripad's axioms (B1) - (B5) characterize the idempotents of any semigroup
relative to the basic products.

Nambooripad noted in his memoir \cite{KSS2} that semilattices are precisely regular biordered sets for which the quasi-orders $\omega^r$ and $\omega^l$ coincide: in this case $\mathcal{S}(e,f) = \{ef\}$. More generally, he characterized regular biordered sets  $E$   for which  $|\mathcal{S}(e,f)| = 1$ for all $e,f \in E$ as {\em pseudo-semilattices} in the sense of Schein \cite{Schein1}. He also provided  characterizations of biordered sets of several other special classes of regular semigroups. In particular, he characterized biordered sets of completely semisimple semigroups, completely regular semigroups, combinatorial regular semigroups, bands, left regular bands, and normal bands.

Several authors have studied biordered sets of various classes of regular semigroups.  Pastijn \cite{Pas2} showed how the biordered set of idempotents of a strongly regular Baer semigroup (in particular the multiplicative semigroup of a regular ring - such as the ring  $M_n(F)$ of $n\times n$ matrices over a field $F$ for example) may be constructed in terms of the complemented modular lattice that is coordinatized by the ring. Pastijn's results play an important role in several papers of Nambooripad that connect regular semigroups with other areas of mathematics (see, for example, the papers \cite{KSS9,NamPas2,KSS7,KrishnaNam}). We will provide a brief description of some of these results in Section 6 below.  In the context of linear algebraic monoids, Putcha \cite{Putcha2} has constructed the biordered set of idempotents of an irreducible linear algebraic monoid from pairs of opposite parabolic subgroups of the group of units of the monoid, which is a reductive group.  There is an extensive literature about biordered sets of locally inverse semigroups: some references to this literature are contained in Section 3 of this paper. Biordered sets of many other classes of regular semigroups have also been studied, see for example the papers \cite{cliff2, Rajan1, Eas2, Eas3, Eas4, EasHall, Pas3,  Prem, KSSKrish}.

\medskip

{\bf Fundamental regular semigroups}

\medskip

A regular semigroup $S$ is said to be {\em fundamental} if the only congruence on $S$ contained in Green's equivalence relation $\mathscr{H}$ is the identity congruence. Munn \cite{Munn1} described fundamental {\em inverse} semigroups with semilattice $E$ as full inverse subsemigroups of the semigroup  of isomorphisms between principal order ideals of $E$. Nambooripad \cite{KSS2} obtained an analogous description of fundamental regular semigroups with biordered set $E$ as a special case of his work on inductive groupoids. Such semigroups are built directly from the biordered set $E$ in much the same way as in Munn's construction, but with some additional complications in the more general case of regular biordered sets. We provide a brief description of his results, as a prelude to his more general theory of inductive groupoids discussed below.

Denote by $T_E$ the set of all $\omega$-isomorphisms of a regular biordered set $E$, that is, all  biorder isomorphisms between sets of the form $\omega(e), \, e \in E$.
If $\alpha :\omega(e)\to\omega(f)$ is an
$\omega $-isomorphism then we write $e_{\alpha}=e$ and $f_{\alpha}=f.$
Clearly  $T_E$  forms a
groupoid under the usual composition of maps: if $\alpha, \beta \in T_E$, then the $\omega$-isomorphism $\alpha\beta$ is defined only if
$f_{\alpha}=e_{\beta}.$

Define the $\omega $-isomorphism $\tau (e,f):\omega(e)\to\omega(f)$
when $e(\mathscr{R} \cup \mathscr{L}) f $
by
\[(g)\tau (e,f)=\begin{cases} gf\text{ if }e\mathscr{R} f\\ fg\text{ if }e\mathscr{L} f \end{cases} \]
for $g\in \omega (e).$
If $\alpha :\omega (e)\to\omega (f)$ is an
$\omega -$isomorphism and if $g\omega e$ then we denote by $g\ast\alpha $ the restriction of  $\alpha$ to
$\omega (g)$.  More generally, if  $g(\omega^r \cup \omega^l) f$ then
\[g\ast \alpha =\begin{cases}\tau(g, ge)(ge\ast\alpha)\text{ if }g\omega^r e\\
\tau (g,eg)(eg\ast\alpha)\text{ if }g\omega^l e. \end{cases}. \]

Similarly for $h(\omega^r\cup\omega^l)f$ then $\alpha\ast h$ is defined by
\[ \alpha\ast h =\begin{cases}(\alpha\ast hf)\tau(hf,h)\text{ if }h\omega^r f\\
(\alpha\ast fh)\tau (fh,h)\text{ if }h\omega^l f. \end{cases}. \]
We regard $\alpha\ast hf$ and $\alpha\ast fh$ as corestrictions in the usual sense.

We can then extend the groupoid
composition in $T_E$ to a semigroup operation. This semigroup operation is
defined on a quotient of $T_E.$ To do this, an equivalence relation $p$ on $T_E$ is defined by
\[ \alpha p\beta\text{ if }e_\alpha\mathscr{R} e_\beta ,\ f_\alpha\mathscr{L} f_\beta\text{ and }
\tau(e_\alpha,e_\beta )\beta =\alpha\tau (f_\alpha,f_\beta).\]

The following theorem of Nambooripad \cite[Theorem 5.2]{KSS2} provides the construction of all fundamental regular
semigroups with a regular biordered set $E$ .
\begin{thm}
\label{fund}
Let $E$ be a regular biordered set and  $p$ be the equivalence  relation on $T_E$
defined above.
Let $[\alpha]$ denote the equivalence class of $\alpha .$
For $\alpha ,\beta\in T_E$ define
\[ [\alpha][\beta]=[(\alpha\ast h)(h\ast \beta)]\]
where $h\in \mathcal{S}(f_\alpha,e_\beta)$.
Then $T_E/p$ is a fundamental regular semigroup whose biordered set  is isomorphic to  $E.$ Conversely, if $S$ is any regular semigroup, then $S$ is fundamental if and only if $S$ is isomorphic to a full regular subsemigroup of  $T_{E(S)}/p$.
\end{thm}

Alternative descriptions of the structure of fundamental regular semigroups were provided by Hall \cite{Hall1} and Grillet \cite{Gri1}. Their approach was extended by Nambooripad \cite{KSS3} in his theory of {\em cross connections}, outlined later in this paper.

\medskip

{\bf Inductive groupoids}

\medskip

Nambooripad made essential use of two  groupoids naturally associated with regular semigroups.  We describe these groupoids below.

Let $S$ be a regular semigroup and let
\[ \mathcal{G}(S)=\{ (x,x'):x'\text{ is an inverse of }x \}.\]
We may view $\mathcal{G}(S)$ as a groupoid with vertex set $v\mathcal{G}(S) = E(S)$  and where $(x,x')$ is considered as a morphism from $xx'$ to $x'x.$ (Here we identify the identity morphism $(xx',xx')$ with the element $xx'$ in $E(S)$.)
The product in the groupoid is defined as follows. For $(x,x'),(y,y')\in \mathcal{G}(S)$
\[ (x,x')(y,y')=(xy,y'x')\text { when }x'x=yy'. \]
The inverse of the morphism $(x,x')$ in $\mathcal{G}(S)$ is of course the morphism $(x',x)$.

There is a natural partial order relation on the groupoid $\mathcal{G}(S)$ defined by
\[(y,y')\le (x,x')\text{ if }yy'\omega xx'\text{ and }y=yy'x,\ y'=x'yy'\]
With respect to this partial order, $\mathcal{G}(S)$ becomes an {\em ordered groupoid} in the sense defined below: ordered groupoids were introduced by Ehresmann \cite{Ehr1} in the context of his work on pseudogroups.

\begin{dfn}\label{og}
Let $\mathcal{G}$ be a groupoid and $\leq$ a partial order on $\mathcal{G}$.  Let $e,f \in v\mathcal{G}$ and $x,y$ denote arbitrary morphisms of $\mathcal{G}$ and let $\mathbf{d}(x)$ and $\mathbf{r}(x)$ denote the domain and codomain respectively of an arbitrary morphism $x$. Then $(\mathcal{G},\leq)$ is called an \emph{ordered groupoid} if the following hold.
\begin{enumerate}
	\item [(OG1)] If $u\leq x$, $v\leq y$ and $\mathbf{r}(u)=\mathbf{d}(v)$, $\mathbf{r}(x)=\mathbf{d}(y)$, then $uv \leq xy$.
	\item [(OG2)] If $x\leq y$, then $x^{-1}\leq y^{-1}$.
	\item [(OG3)] If $1_e\leq 1_{\mathbf{d}(x)}$, then there exists a unique element $e{\downharpoonleft} x$ (called the \emph{restriction} of $x$ to $e$) in $\mathcal{G}$ such that $e{\downharpoonleft} x\leq x$ and $\mathbf{d}(e{\downharpoonleft} x) = e$. 		
	\item[(OG3$^*$)] If $1_f\leq 1_{\mathbf{r}(x)}$, then there exists a unique element $x{\downharpoonright} f$ (called the \emph{corestriction} of $x$ to $f$) in $\mathcal{G}$ such that $x{\downharpoonright} f \leq x$ and $\mathbf{r}(x{\downharpoonright} f) = f$.
\end{enumerate}
\end{dfn}
In fact axioms (OG1), (OG2) and (OG3) are  equivalent to (OG1), (OG2) and (OG3$^*$).

An {\em $E$-sequence} in a regular biordered set $E$ is  a  sequence $e_1,e_2,...,e_n$ of elements of  $E$ such that $e_k (\mathscr{R} \cup \mathscr{L}) e_{k+1}$ for $k = 1,...,n-1$. An element $e_i$ of such a sequence is called {\em inessential} if either $e_{i-1} \mathscr{R} e_i \mathscr{R} e_{i+1}$ or $e_{i-1} \mathscr{L} e_i \mathscr{L} e_{i+1}$. The unique $E$-sequence $[e_1,e_2,...,e_n]$ obtained by removing all inessential vertices from an $E$-sequence $(e_1,e_2,...,e_n)$ is called an {\em $E$-chain}. Denote the set of all $E$-chains in $E$ by $\mathcal{G}(E)$. Then $\mathcal{G}(E)$ is a groupoid with vertex set $v\mathcal{G}(E) = E$ with respect to the multiplication

\[[e_1,e_2,...,e_n][f_1,f_2,...,f_m] = [e_1,e_2,...e_n,f_2,...,f_m],\]

\noindent defined if and only if $e_n = f_1$. Here ${\bf d}([e_1,e_2,...e_n]) = e_1, \, {\bf r}([e_1,e_2,...,e_n]) = e_n$ and the inverse of  $[e_1,...,e_n]$ is $[e_n,...,e_1]$.

If $c=[e_1,e_2,...,e_n]$ is an $E$-chain and $h \in \omega^r(e_1)$ then we define $h*c$ to be the $E$-chain
\[h*c =  [h=h_0,h_1,h_2,...,h_n]\text{ where }h_i = e_ih_{i-1}e_i,\, i = 1,...,n\]
Dually, if $k \in \omega^l(e_n)$ then we define
\[c*k = [k_1,k_2,...,k_n,k_{n+1}=k]\text{ where }k_i = e_ik_{i+1}e_i, \, i = i,...,n\]
If $c = [e_1,...,e_n]$ and $c' = [f_1,...f_m]$ are two $E$-chains, then we define
\[c' \leq c\text{ if and only if }f_1 \omega e_1\text{ and }c'= f_1*c\]
Then, with respect to this partial order, $\mathcal{G}(E)$ is an ordered groupoid, called the {\em groupoid of $E$-chains of $E$}.

Now suppose that $S$ is a regular semigroup with biordered set $E(S)$ and that $c = [e_1,e_2,...,e_n]$ is an $E$-chain in $E(S)$. Then it follows from the Miller-Clifford theory of regular $\mathscr{D}$-classes  that the (non-basic) product $e_1e_2...e_n$ is in the $\mathscr{H}$-class $R_{e_1} \cap L_{e_n}$ and similarly that $e_n...e_2e_1 \in L_{e_1} \cap R_{e_n}$. In fact $e_n...e_2e_1$ is an inverse of $e_1e_2...e_n$. It follows that   $\epsilon_{S}(c) = (e_1e_2...e_n,e_n...e_2e_1) \in \mathcal{G}(S)$. We regard $ \epsilon_{S}$ as a functor from $\mathcal{G}(E)$ to $\mathcal{G}(S)$ that enables us to {\em evaluate} the $E$-chain $c$ in the groupoid $\mathcal{G}(S)$.

Nambooripad identified one more essential concept that is naturally associated with a biordered set and that is needed in his definition of an inductive groupoid. Given a biordered set $E$, a $2 \times 2$ matrix $\bigl[ \begin{smallmatrix} e&f\\ g&h \end{smallmatrix} \bigr]$ of elements of $E$ such that $e\mathrel{\mathscr{R}} f\mathrel{\mathscr{L}}h\mathrel{\mathscr{R}} g\mathrel{\mathscr{L}}e$ forms a distinguished E-chain and is known as an \emph{E-square}. An E-square of the form  $\bigl[ \begin{smallmatrix} g&h\\ eg&eh \end{smallmatrix} \bigr]$
where $g,h \in \omega^l(e)$ and $g\mathrel{\mathscr{R}} h$ is said to be row-singular. Dually, we define column singular E-squares and an E-square is said to be \emph{singular} if it is either row-singular or column-singular. It is not difficult to see (once it is pointed out) that if $\bigl[ \begin{smallmatrix} e&f\\ g&h \end{smallmatrix} \bigr]$ is a singular $E$-square in the biordered set $E(S)$ of a  semigroup $S$ and $\epsilon = \epsilon_{S}: \mathcal{G}(E) \to \mathcal{G}(S)$ is the functor defined above, then the $E$-square $\bigl[ \begin{smallmatrix} e&f\\ g&h \end{smallmatrix} \bigr]$ is {\em $\epsilon$-commutative}, that is

$$\epsilon(e,f)\epsilon(f,h) =\epsilon(e,g)\epsilon(g,h).$$

The preceding discussion motivates Nambooripad's concept of an {\em inductive groupoid} defined below. A functor $F$ between two ordered groupoids is said to be a $v$-isomorphism if the object map $vF$ is an order isomorphism.

\begin{dfn}\label{dfnig}
Let $E$ be a regular biordered set, $\mathcal{G}$ an ordered groupoid and let $\epsilon\colon \mathcal{G}(E)\to \mathcal{G}$ be a $v$-isomorphism called an \emph{evaluation functor}. We say that $(\mathcal{G},\epsilon)$ forms an \emph{inductive groupoid} if the following axioms and their duals hold.
\begin{enumerate}
	\item[(IG1)] Let $x\in \mathcal{G}$ and for $i=1,2$, let $e_i$, $f_i \in E$ such that $\epsilon (e_i) \leq \mathbf{d}(x)$ and $\epsilon (f_i) = \mathbf{r}(\epsilon (e_i){\downharpoonleft} x)$. If $e_1\mathrel{\omega}^r e_2$, then $f_1\mathrel{\omega}^r f_2$, and
	$$\epsilon (e_1,e_1e_2)(\epsilon (e_1e_2){\downharpoonleft} x) = (\epsilon (e_1){\downharpoonleft} x)\epsilon (f_1,f_1f_2).$$
	\item[(IG2)] All singular E-squares are $\epsilon $-commutative.
\end{enumerate}
\end{dfn}

Nambooripad \cite{KSS2} defines morphisms between inductive groupoids as follows. Let $(\mathcal{G},\epsilon )$ and $(\mathcal{G}',\epsilon ')$ be two inductive groupoids with biordered sets $E$ and $E'$ respectively. An order preserving functor $F\colon \mathcal{G} \to\mathcal{G}'$ is said to be \emph{inductive} if $vF\colon E \to E'$ is a regular bimorphism of biordered sets such that the following diagram commutes.
\begin{equation}\label{indf}
\xymatrixcolsep{4pc}\xymatrixrowsep{3pc}\xymatrix
{
\mathcal{G}(E) \ar[r]^{\mathcal{G}(vF)} \ar[d]_{\epsilon }
& \mathcal{G}(E') \ar[d]^{\epsilon '} \\
\mathcal{G} \ar[r]^{F} & \mathcal{G}'
}
\end{equation}
Then the class of inductive groupoids with inductive functors as morphisms forms a category.

If $S$ is a regular semigroup, then the pair $(\mathcal{G}(S),\epsilon_{S})$ is an inductive groupoid. Conversely, suppose that $\mathcal{G} = (\mathcal{G},\epsilon)$ is an inductive groupoid. We now outline Nambooripad's construction of an associated regular semigroup $S(\mathcal{G})$.

The evaluation functor enables one to extend the restrictions
$e{\downharpoonleft} x$ defined in ordered groupoids to cases where $e\omega^r {\bf d}(x)$ and $e\omega^l{\bf d}(x)$
These are defined as follows.
\[e\ast x=\epsilon (e,e{\bf d}(x))(e{\bf d}(x){\downharpoonleft} x)\text{ if }e\omega^r {\bf d}(x)\]
and
\[e\ast x=\epsilon (e,{\bf d}(x)e)({\bf d}(x)e{\downharpoonleft} x)\text{ if }e\omega^l {\bf d}(x).\]
Similarly corestrictions can also be extended. The semigroup associated with an inductive groupoid is a quotient
$\mathcal{G}/p$ where $p$ is an equivalence relation on $\mathcal{G}$ defined by
\[x p y\text{ if and only if }{\bf d}(x) \mathscr{R} {\bf d}(y), \, {\bf r}(x) \mathscr{L} {\bf r}(y)\text{ and }x\epsilon ({\bf r}(x),{\bf r}(y)) = \epsilon ({\bf d}(x),{\bf d}(y))y.\]

Nambooripad \cite[Theorems 4.12 and 4.14]{KSS2} proves the following  structure theorem for regular semigroups in terms of inductive groupoids.

\begin{thm}\label{thmind}
\label{indreg}
Let $(\mathcal{G},\epsilon)$ be an inductive groupoid and let $p$ be the equivalence defined above.
Let $[x]$ denote the equivalence class of $x.$ For $x,y\in \mathcal{G}$ define
\[ [x][y]=[(x\ast h)(h\ast y)]\]
where $h\in \mathcal{S}({\bf r}(x),{\bf d}(y))$.
Then this binary operation is well defined and $S(\mathcal{G}) = \mathcal{G}/p$ is a regular semigroup with this binary operation.
Further the biordered set of $S(\mathcal{G})/p$ is isomorphic to the biordered set of $\mathcal{G}$. Conversely, if $S$ is any regular semigroup, then $S$ is isomorphic to $S(\mathcal{G}(S))$. The category of regular semigroups is equivalent to the category of inductive groupoids.
\end{thm}

\begin{rmk}
If the biordered set $E$ is a semilattice (i.e. $\omega^r = \omega^l$), then most of the complication required to define inductive groupoids vanishes. In this case $\mathcal{G}(E)$ is the trivial groupoid
consisting of identity morphisms alone, $\epsilon_{\mathcal G}$ is just the identity map on $E$
and there are no non-trivial $E$-squares, so in this case every ordered groupoid is inductive. Furthermore, $\mathcal{S}({\bf r}(x),{\bf d}(y)) = \{x^{-1}xyy^{-1}\}$ so Theorem \ref{indreg} reduces to Schein's construction \cite{Schein2} of inverse semigroups via inductive groupoids. This result was referred to as the ``Ehresmann-Schein-Nambooripad"(ESN) theorem  by Lawson in his book \cite{lawson}. Furthermore, in this case, the element $h*y$ in the definition of the multiplication in Theorem \ref{indreg} is the image of $y$ under the {\em structure mapping} from $R_y$ to $R_h$ (and $x*h$ is  the image of $x$ under the structure mapping from $L_{x}$ to $L_h$), so Theorem \ref{indreg} also reduces to Meakin's equivalent construction \cite{Me1} of inverse semigroups via structure mappings.

\end{rmk}

\begin{rmk}
It is evident that the construction of fundamental regular semigroups given in Theorem \ref{fund} may be viewed as a special case of Theorem \ref{indreg}. The corresponding evaluation functor is the natural extension to $E$-chains of the $\omega$-isomorphisms $\tau(e,f)$ for $e (\mathscr{R} \cup \mathscr{L}) f$. That is, $\tau([e_1,e_2,e_3,...,e_{n-1},e_n]) = \tau(e_1,e_2)\tau(e_2,e_3)...\tau(e_{n-1},e_n)$.
\end{rmk}


\section{Outgrowths of Nambooripad's inductive groupoid approach}\label{secindog}

In this section, we briefly discuss several concepts and outgrowths of Nambooripad's inductive groupoid approach to the structure of regular semigroups.

\medskip
{\bf The natural partial order on a regular semigroup}

\medskip
In his paper \cite{KSS4}, Nambooripad defined a natural partial order on a regular semigroup as follows. Let $S$ be a regular semigroup and $x,y \in S$. Then define
\[x \leq y\text{ if and only if }R_x \leq R_y\text{ and }x = fy\text{ for some idempotent }f \in R_x.\]

It is not too difficult to show (see \cite[Proposition 1.2]{KSS4}) that the relation $\leq$ is self-dual and that this relation is a partial order on $S$ that coincides with the well-known definition of the natural partial order on an inverse semigroup, if $S$ is inverse. Furthermore, if $e\in E(S)$ and $f \leq e$ then $f \in E(S)$ and the restriction of the relation $\leq$ to $E(S)$ is the relation $\omega$. One observes that if $x*h$ and $h*y$ are the elements that appear in the multiplication in the statement of Theorem \ref{indreg}, then $(x*h) \leq x$ and $(h*y) \leq y$. The element $x*h$ is the image of $x$ under the structure mapping from $L_x$ to $L_h$  in the sense of Meakin \cite{Me2} (and a dual statement applies to $h*y$).

Nambooripad shows in \cite{KSS4} that the natural partial order $\leq$ on a regular semigroup $S$ enjoys many of the same properties as the natural partial order on an inverse semigroup. However, unlike the situation for inverse semigroups, the relation $\leq$ on a regular semigroup is not in general compatible with the multiplication in $S$. In fact Nambooripad shows that this is the case if and only if $S$ is a {\em pseudo-inverse} semigroup, that is $E(S)$ is a pseudo-semilattice in the sense of Schein \cite{Schein1}.  When combined with several of his other results from his papers \cite{KSS2, KSS5, KSS6}, Nambooripad proves the following theorem.

\begin{thm}
\label{pseudoinv}

The following conditions on a regular semigroup $S$ are equivalent.

(a) $S$ is a pseudo-inverse semigroup (i.e. $E(S)$ is a pseudo-semilattice).

(b) $S$ is a locally inverse semigroup (i.e. $eSe$ is an inverse semigroup for each $e \in E(S)$).

(c) $\omega(e)$ is a semilattice for each $e \in E(S)$.

(d) $|\mathcal{S}(e,f)| = 1$ for all $e,f \in E(S)$.

(e) The natural partial order $\leq$ on $S$ is compatible with the multiplication in $S$.

(f) If $x \leq y$ in $S$ then for every $(y_1,y_2) \in L_y \times R_y$ there is a unique pair $(x_1,x_2) \in L_x \times R_x$ such that $x_i \leq y_i$ for $i = 1,2$.

\end{thm}

There is a large literature devoted to the study of locally inverse (pseudo-inverse) semigroups and pseudosemilattices (also called local semilattices by some authors). In addition to Nambooripad's papers \cite{KSS2,KSS4,KSS5,KSS6} we refer the reader to the papers \cite{Au, Au2,  AuOl, BillSz,  BMP, Kad, McA1, Me3, Me4,  MeKSS1, MePas1, MePas2, Oli, Pas1, PasPet, Sze1, Veer } for much additional information about the structure of pseudosemilattices and locally inverse semigroups.

\medskip

{\bf Regular idempotent-generated semigroups}

\medskip

Idempotent-generated semigroups have been the subject of much study in the literature. For example, an early result of J.A. Erd\"os \cite{Erd} shows that the idempotent-generated part of the semigroup $M_n(F)$ of $n \times n$ matrices over a field $F$ consists of the identity matrix and all singular matrices. J.M. Howie \cite{Howie} proved a similar result for the full transformation monoid on a finite set and also showed that every semigroup may be embedded in a suitable idempotent-generated semigroup. These results have been extended in many different ways and many authors have studied the structure of idempotent-generated semigroups. For example, in a significant extension of Erd\"os' result,  Putcha \cite{Putcha1} gave necessary and sufficient conditions for a reductive linear algebraic monoid to have the property that every non-unit is a product of idempotents.

In his memoir \cite{KSS2}, Nambooripad provided a construction of all regular idempotent-generated semigroups with a given biordered set.
He first observed that inductive groupoids with surjective evaluations completely characterize regular idempotent-generated semigroups. This follows  from the fact that if $x$ is a regular element in the idempotent-generated part of a semigroup $S$, then there is an $E$-chain $[e_1,e_2,...,e_n]$ such that $x = e_1e_2...e_n$. (This was observed by Fitz-Gerald \cite{FitzG} and also by Nambooripad as a consequence of Theorem 4.13 of his memoir \cite{KSS2}). Nambooripad then studied properties of {\em $E$-cycles} to construct all regular idempotent-generated semigroups with a given biordered set $E$. We briefly summarize his ideas.

An $E$-cycle at $e$ is an $E$-chain $c = [e_1,e_2,...,e_n]$ with $e = e_1 = e_n$. Clearly every $E$-square $\bigl[ \begin{smallmatrix} e&f\\ g&h \end{smallmatrix} \bigr]$ gives rise to an $E$-cycle $[e,f,h,g,e]$ at $e$: a {\em singular $E$-cycle} is an $E$-cycle determined by a singular square. We denote the set of singular $E$-cycles by $\Gamma_0$. Let $E$ be the underlying biordered set of an inductive groupoid $\mathcal{G} = (\mathcal{G},\epsilon)$. We say that an $E$-cycle $c = [e,e_2,...e_{n-1},e]$ is {\em $\epsilon$-commutative} if $\epsilon(c) = \epsilon(e,e)$. Denote the set of all $\epsilon$-commutative $E$-cycles by $\Gamma_{\epsilon}$. As noted in Section 2,  every singular $E$-cycle $c$ at $e$ is $\epsilon$-commutative.   It also follows from \cite{KSS2} that every $\epsilon$-commutative $E$-cycle is $\tau$-commutative; that is, $\Gamma_0 \subseteq \Gamma_{\epsilon} \subseteq \Gamma_{\tau}$. In fact $\Gamma_{\epsilon}$ is a {\em proper} set of $E$-cycles in the sense of the following definition.

\begin{dfn}
\label{properdef}
A set $\Gamma$ of $E$-cycles of a biordered set $E$ is said to be proper if it satisfies the following:
\begin{enumerate}
	\item[(P1)] $\Gamma_0 \subseteq \Gamma \subseteq \Gamma_{\tau}$;
	\item[(P2)] $\gamma \in \Gamma$ implies $\gamma^{-1} \in \Gamma$;
	\item[(P3)] If $\gamma$ is an $E$-cycle at $e$ in $\Gamma$ and $f \in \omega(e)$ then $f*\gamma \in \Gamma$.
\end{enumerate}
\end{dfn}

If $(\mathcal{G},\epsilon)$ is an inductive groupoid  then it is clear that if $c_1 = [e_1,e_2,...,e_n]$ and $c_2 = [f_1,f_2,...,f_m]$ are $E$-chains with $f_1 = e_n$ and if $\gamma=[g_1,g_2,...,g_r]$ is an $\epsilon$-commutative $E$-cycle at $e_n (= f_1)$, then $\epsilon (c_1c_2) = \epsilon (c_1 \gamma c_2)$: that is, these two $E$-chains have the same $\epsilon$-value in $S(\mathcal{G})$ and hence the products $w(c_1c_2) = e_1e_2...e_nf_1f_2...f_m$ and $w(c_1 {\gamma} c_2) = e_1e_2...e_ng_1g_2...g_rf_1f_2...f_m$ coincide in $S(\mathcal{G})$. We write this as $c_1 c_2 \mapsto_{\Gamma} c_1 \gamma c_2$. In general, for any proper set $\Gamma$ of $E$-cycles, and $c,c'$ any two $E$-chains, we write $c \sim_{\Gamma} c'$ if there is some sequence of $E$-chains $c = c_1,c_2,...,c_n=c'$ such that  $c_i \mapsto_{\Gamma} c_{i+1}$  for each $i = 1,...,n-1$. Then $\sim_{\Gamma}$ is an equivalence relation and  $w(c) = w(c')$ in $S(\mathcal{G})$ if $c \sim_{{\Gamma}} c'$. Let $\mathcal{G}_{\Gamma} = \mathcal{G}(E)/{\sim_{\Gamma}}$ and let $\epsilon_{\Gamma}$ be the canonical surjection of $\mathcal{G}(E)$ onto $\mathcal{G}_{\Gamma}$. Nambooripad proves the following result based on Corollary 6.8 and Theorem 6.9 of his memoir \cite{KSS2}.

\begin{thm}
\label{rig}
Let $\Gamma$ be a proper set of $E$-cycles in a regular biordered set $E$. Then $\mathcal{G}_{\Gamma}$ is an inductive groupoid with surjective evaluation $\epsilon_{\Gamma}$ so $B_{\Gamma} = S(\mathcal{G}_{\Gamma})$ is a regular idempotent-generated semigroup with biordered set $E$. Every regular idempotent-generated semigroup with biordered set $E$ can be constructed this way. If $\mathcal{G}$ is an inductive groupoid with surjective evaluation $\epsilon$, then $\mathcal{G}$ is isomorphic to $\mathcal{G}_{\Gamma_{\epsilon}}$.

\end{thm}

Nambooripad provides a somewhat more detailed construction of the semigroups $B_{\Gamma}$ and  more information about their properties in Theorems 6.9 and 6.10 of his memoir \cite{KSS2}. It is a consequence of his results that the semigroup $RIG(E) = B_{\Gamma_0}$ corresponding to the  set $\Gamma_0$ of singular squares of a regular biordered set $E$ is a universal object in the category of regular idempotent-generated semigroups with biordered set $E$. The semigroup $RIG(E)$ has been much studied in subsequent literature, where it is referred to  as the {\em free regular idempotent-generated semigroup} on the biordered set $E$. It is a homomorphic image of the {\em free idempotent-generated semigroup $IG(E)$} on $E$ which is the semigroup with presentation

\[IG(E) = \langle E : e \cdot f = ef\text{ if } ef\text{ is a basic product in }E \rangle\]

Here $e \cdot f$ is a word of length $2$ in $E^*$ and $ef$ is a word of length $1$. The semigroup $RIG(E)$ is obtained from $IG(E)$ by adding the relations

\[e \cdot f = e \cdot h \cdot f\text{ for all }e,f \in E\text{ and }h \text{ in }\mathcal{S}(e,f).\]

While the semigroup $RIG(E)$ corresponding to a regular biordered set is in a sense known from Nambooripad's construction of its inductive groupoid, the structure of this semigroup is not well  understood. In particular, the study of the structure of the maximal subgroups of $RIG(E)$ has received considerable attention in the literature. It is known (see \cite{BMM1}) that the maximal subgroup of $RIG(E)$ containing an idempotent $e \in E$ is isomorphic to the maximal subgroup of $IG(E)$ containing $e$. In their paper \cite{NamPas1}, Nambooripad and Pastijn showed that if $E$ is a biordered set that has no non-degenerate singular squares (for example a locally inverse semigroup), then the maximal subgroups of $RIG(E)$ are free groups. In fact it was expected that the maximal subgroups of $RIG(E)$ should be free groups for any regular biordered set $E$: this was explicitly conjectured  in a paper  by McElwee \cite{McEl}. However, in their paper \cite{BMM1}, Brittenham, Margolis and Meakin used topological methods to construct an example of a finite regular biordered set $E$ for which the maximal subgroups of $RIG(E)$ are free abelian groups of rank $2$. Subsequently it was shown by Gray and Ruskuc \cite{GrayRus1}, using Reidemeister-Schreier rewriting methods that {\em every} group arises as a maximal subgroup of $RIG(E)$ for some biordered set $E$. An elegant alternative proof of this was provided by Gould and Yang \cite{GouldYang1}.

It was shown in \cite{BMM2} that if $E$ is the biordered set of the monoid $M_n(F)$ of $n \times n$ matrices over a field (or division ring) $F$, then the maximal subgroup of $RIG(E)$ corresponding to an idempotent matrix of rank $1$ is isomorphic to the multiplicative group of units of $F$ and the maximal subgroup of $M_n(F)$ corresponding to an idempotent matrix of rank $n-1$ is a free group. These results were extended significantly by Dolinka and Gray \cite{DG1} who showed that the maximal subgroup of $RIG(E)$ corresponding to an idempotent matrix of rank $r$ is isomorphic to $GL_r(F)$ provided $1 \leq r < n/3$. However the structure of these groups corresponding to idempotent matrices of rank $r$ with $n/3 \leq r < n-1$ remains unknown. Further results on maximal subgroups of free  idempotent-generated semigroups over various biordered sets may be found in several papers, for example \cite{Do1, Do2, DGY, DoR, GrayRus2}.

Recent work in this area has been concerned with the study of the word problem and additional structure of the semigroups $IG(E)$ and $RIG(E)$. The reader is referred to the papers \cite{DGY2, GY2, Do3, DGR  } for much work along these lines. It seems that little is known about the structure of  Nambooripad's semigroups $B_{\Gamma}$ where $\Gamma$ is a proper set of $E$-cycles other than $\Gamma_0$ or $\Gamma_{\tau}$.



\medskip

{\bf Beyond regular semigroups}

\medskip

Additional literature connected to Nambooripad's work on the inductive groupoid construction
of regular semigroups is concerned with
various generalizations of the ``ESN" theorem.
Different non-regular generalisations of regular and inverse semigroups including concordant semigroups, abundant semigroups, Ehresmann semigroups, ample semigroups, restriction semigroups and weakly $U$-regular semigroups have been described using generalisations of inductive groupoids, including the concepts of  inductive cancellative categories, Ehresmann categories, inductive categories, inductive constellations, and weakly regular categories. We refer the reader to the papers  \cite{armstrong,lawsonordered0,lawsonenlarge,lawsonordered,gomes,gouldrestr,hollings2,gouldrestriction,gould2,hollings,wang,wangureg} for much work along these lines.

\section{Regular semigroups and cross connections} \label{seccxn}

As mentioned earlier, Grillet \cite{Gri1} developed the ideas initiated by Hall \cite{Hall1} exploiting the ideal structure of a regular semigroup to give a construction of fundamental regular semigroups via the notion of cross-connections. Recall that the principal left and right ideals of a regular semigroup $S$ are given by the following sets $\mathcal{L}(S)$ and $\mathcal{R}(S)$, respectively:

$$\mathcal{L}(S)=\{Se: e\in E(S)\} \text{ and } \mathcal{R}(S)=\{eS: e\in E(S)\}.$$

Under the usual set inclusion, these sets form partially ordered sets (posets). Grillet characterised these posets as {\em regular posets} and showed that $\mathcal{L}(S)$ and $\mathcal{R}(S)$ are inter related via a {\em fundamental cross-connection}. Conversely, given a pair of abstractly defined regular posets and a fundamental cross-connection, one could obtain a fundamental regular semigroup.

In 1978, Nambooripad \cite{bicxn} proved that  his construction of fundamental regular semigroups via biordered sets (as described in Section \ref{secind}) and Grillet's construction via cross-connections (as briefly indicated above) are equivalent. Recall that Nambooripad's construction of fundamental regular semigroups may be realised as a special case of his inductive groupoid construction. In other words, by `attaching' a groupoid structure to the biordered set, one could construct arbitrary regular semigroups. So, it was naturally enticing to `attach' a category structure to cross-connections and try to construct arbitrary regular semigroups. Although such an attempt necessitated the development of a rather complicated machinery, Nambooripad went forward with this plan. A major reason behind this effort is the fact that such an approach will overcome the barrier imposed by `idempotents' and lead towards a much more general framework for structure theorems. This aspect of cross-connections shall be elaborated on  in the next section.

In this section, we proceed to describe how Nambooripad \cite{cross0,KSS3} generalised  Grillet's cross-connections by employing {\em normal categories}. Normal categories are categorical abstractions of the set of  principal one-sided ideals of a regular semigroup. First, we shall see how a regular semigroup gives rise to a cross-connection between its left and right normal categories. Conversely, we shall briefly describe how an abstractly defined cross-connection between two normal categories will give rise to a regular semigroup.

\medskip

{\bf Normal categories}

\medskip

We begin by describing the normal category $\mathbb{L}(S)$ of principal left ideals of a regular semigroup $S$: this will lead to its axiomatisation; dually we can characterise $\mathbb{R}(S)$.

Given regular semigroup $S$, the category $\mathbb{L}(S)$ that arises from the principal left ideals is given by:
$$v\mathbb{L}(S) = \{ Se : e \in E(S) \},$$
and for each $x\in Se$ and for each $u\in eSf$, a morphism from $Se$ to $Sf$ is the function $\rho(e,u,f)\colon x \mapsto xu $. Thus the set of all morphisms in  the category $\mathbb{L}(S)$ from the object $Se$ to $Sf$ is given by the set
$$ \mathbb{L}(S)(Se,Sf) = \{ \rho(e,u,f) : u\in eSf \}.$$
Given any two morphisms, say $\rho(e,u,f)$ and $\rho(g,v,h)$, they are equal if and only if $e\mathrel{\mathscr{L}}g$, $f\mathrel{\mathscr{L}}h$ and $v=gu$. Two morphisms $\rho(e,u,f)$ and $\rho(g,v,h)$ are composable if $Sf=Sg$ (i.e., if $f\mathrel{\mathscr{L}}g$) and then
$$\rho(e,u,f)\:\rho(g,v,h) = \rho(e,uv,h).$$

Observe that the set $v\mathbb{L}(S)$ is exactly the same (regular) poset explored by Grillet and this poset can be identified as a subcategory of $\mathbb{L}(S)$, via the distinguished morphisms of the form $\rho(e,e,f)$, whenever $Se\subseteq Sf$. These morphisms correspond to inclusion maps and hence are called {\em inclusions}. Nambooripad used the notion of \emph{category with subobjects} to abstract such a category, wherein a poset sits inside as a distinguished $v$-full subcategory. In a category $\mathcal{C}$ with subobjects, if $c,d \in v\mathcal{C}$ such that $c\le d$, then there is an inclusion from $c \to d$ and in the sequel, we denote this inclusion morphism by $j(c,d)$.

For an inclusion $\rho(e,e,f)\in\mathbb{L}(S)$, it is easy to see that $\rho(f,fe,e)$ is a right inverse; then we say that the inclusion \emph{splits} and the inverse morphism $\rho(f,fe,e)$ is called a {\em retraction}.

Further observe that any morphism $\rho(e,u,f)$ in $\mathbb{L}(S)$ can be factorised as
$$\rho(e,u,f)= \rho(e,g,g)\rho(g,u,h)\rho(h,h,f)$$
for some $h\in E(L_u)$ and $g\in E(R_u)\cap \omega(e)$, so that $\rho(e,g,g)$ is a retraction, $\rho(g,u,h)$ is an isomorphism and $\rho(h,h,f)$ is an inclusion. Hence, given a morphism $f$ in a category $\mathcal{C}$ with subobjects, a factorisation of the form $f=quj$ is called as a \emph{normal factorisation} if $q$ is a retraction, $u$ is an isomorphism and $j$ is an inclusion. Here, the epimorphism $qu$ is called the epimorphic part of $f$ and is denoted by $f^{\circ}$. The codomain of $f^{\circ}$ is known as the image of $f$ and is denoted by im $f$. Indeed, it is this factorisation property which replaces the role of restriction/corestriction of the inductive groupoid construction.


Recall that our aim is to develop a framework to build the semigroup back from the abstractly defined structures. For this, Grillet used the certain mappings on the regular posets called {\em normal mappings}. The basic building block of Nambooripad's construction was provided using the following notion of a {\em normal cone} which may be viewed as an extension of a {normal mapping}.

\begin{dfn}
	Let $\mathcal{C}$ be a category with subobjects and $d\in v\mathcal{C}$. Then a function $\gamma\colon v\mathcal{C}\to \mathcal{C}$, $a\mapsto\gamma(a)\in\mathcal{C}(a,d)$ is said to be a \emph{normal cone} with vertex $d$ if:
	\begin{enumerate}
		\item whenever $a\subseteq b$, $j(a,b) \gamma(b)=\gamma(a)$;
		\item there exists at least one $c\in v\mathcal{C}$ such that $\gamma(c)\in\mathcal{C}(c,d)$ is an isomorphism.
	\end{enumerate}
\end{dfn}

Essentially, normal cones are certain `pastings' of the morphisms but they are direct abstractions of the right regular representation of a semigroup. For instance, if $a$ is an arbitrary element of $S$, then for each $Se\in v\mathbb{L}(S)$,  the function $\rho^a\colon v\mathbb{L}(S)\to \mathbb{L}(S)$ defined by
\begin{equation}\label{eqnprinc}
\rho^a(Se) = \rho(e,ea,f)\text{ where } f\in E(L_a)
\end{equation}
is a normal cone with vertex $Sf$, usually referred to as a {\em principal cone}. Observe that, for an idempotent $e\in E(S)$, we have a principal cone $\rho^e$ with vertex $Se$ such that $\rho^e(Se)= \rho(e,e,e)=1_{Se}$. In fact, it can be shown that the collection $\{\rho^a: a\in S\}$ of all principal cones in $\mathbb{L}(S)$ is isomorphic to the right regular representation $S_\rho$ of the semigroup $S$.


Summarising the above discussion leads us to the following abstraction of the category $\mathbb{L}(S)$ of principal left ideals of a regular semigroup $S$.
\begin{dfn}
	\label{dfnnormc}
	A small category $\mathcal{C}$ is said to be a \emph{normal category} if:
	\begin{enumerate}
		\item [(NC 1)] $\mathcal{C}$ is a category with subobjects;
		\item [(NC 2)] every inclusion in $\mathcal{C}$ splits;
		\item [(NC 3)] every morphism in $\mathcal{C}$ admits a normal factorisation;
		\item [(NC 4)] for each $c\in v\mathcal{C}$ there exists a normal cone $\mu$ such that $\mu(c)=1_c$.
	\end{enumerate}
\end{dfn}

Now, given such an abstractly defined normal category $\mathcal{C}$, Nambooripad introduces a special binary operation on the set of all normal cones in $\mathcal{C}$. Observe that for a normal cone $\gamma$ with vertex $c_\gamma$ and an epimorphism $f\in\mathcal{C}(c_\gamma,d)$, we can construct a new normal cone $\gamma\ast f$ with vertex $d$ such that
\[\gamma\ast f\colon a\mapsto\gamma (a)f \]
for each $a\in v\mathcal{C}$. So, given any two normal cones $\gamma ,\delta$ in $\mathcal{C}$, we define the product
\[\gamma\:\delta =\gamma \ast [\delta (c_\gamma )]^{\circ}\]
where $[\delta (c_\gamma )]^{\circ}$ is the epimorphic part of the morphism $\delta (c_\gamma )$. Then the set of all normal cones in $\mathcal{C}$, denoted in the sequel by $T ({\mathcal{C}})$, forms a regular semigroup. Further, the normal category $\mathbb{L}(T ({\mathcal{C}}))$ of the principal left ideals of the regular semigroup $T ({\mathcal{C}})$, is isomorphic to $\mathcal{C}$. Hence we have:
\begin{thm}\cite[Theorems III.16 and III.19]{KSS3}\label{thmls}
	A small category $\mathcal{C}$ is normal if and only if $\mathcal{C}$ is isomorphic to the category $\mathbb{L}(S)$, for some regular semigroup $S$.
\end{thm}

Dually, we define the normal category $\mathbb{R}(S)$ of principal right ideals of a semigroup by:
$$v\mathbb{R}(S) = \{ eS : e \in E(S) \} \: \text{ and }  \: \mathbb{L}(S)(eS,fS) = \{ \lambda(e,u,f) : u\in fSe \}$$
where for each $x\in eS$ and for each $u\in fSe$, a morphism from $eS$ to $fS$ is the function $\lambda(e,u,f)\colon x \mapsto ux$.

Having characterised the principal ideal structure of a regular semigroup as the normal categories $\mathbb{L}(S)$ and $\mathbb{R}(S)$, we proceed to investigate their inter-relationship. This rather non-trivial connection is captured using the notion of a cross-connection and that requires the introduction of a dual category associated with a given normal category.

\medskip

{\bf Normal duals and cross-connections}

\medskip

In \cite{Gri1}, Grillet  devised a dual of a regular poset by using a certain set of equivalence relations on the given poset. To extend Grillet's construction to categories, Nambooripad defined certain special set valued functors called $H$-functors. Given any small category $\mathcal{C}$, it is well known that there is an associated functor category $[\mathcal{C},\mathbf{Set}] $ with the set of objects as all functors from $\mathcal{C}$ to the category $\mathbf{Set}$ and  natural transformations as morphisms.

Let $\gamma$ be a normal cone in $T (\mathcal{C})$. Then for each $c\in v\mathcal{C}$ and $g\in\mathcal{C}(c,d)$, Nambooripad defined an $H$-functor $H(\gamma;-)\colon \mathcal{C}\to \mathbf{Set}$ as follows:
\begin{equation} \label{eqnH}
\begin{split}
H({\gamma};{c})&= \{\gamma\ast f^\circ : f \in \mathcal{C}(c_{\gamma},c)\} \text{ and }\\
H({\gamma};{g}) \colon H({\gamma};{c}) &\to H({\gamma};{d}) \text{ given by }\gamma\ast f^\circ \mapsto \gamma\ast (fg)^\circ
\end{split}
\end{equation}
It was shown that $H(\gamma;-)$ is a representable functor such that there exists an associated natural isomorphism $\eta_\gamma\colon H(\gamma;-) \to \mathcal{C}(c_\gamma,-)$ where $\mathcal{C}(c_\gamma,-)$ is the covariant hom-functor determined by the object $c_\gamma$.

Thus, given a normal category $\mathcal{C}$, Nambooripad defined the \emph{normal dual} N$^*\mathcal{C}$  as the full subcategory of $[\mathcal{C},\mathbf{Set}] $ such that
$$vN^*\mathcal{C}=\{H(\gamma;-):\gamma\in T (\mathcal{C})\}.$$
So the morphisms in the normal dual N$^*\mathcal{C}$ are natural transformations between the $H$-functors. Using this machinery, Nambooripad proved the following \cite[Theorem III.25]{KSS3}:
\begin{thm}\label{thmdualrs}
Let $\mathcal{C}$ be a normal category with the normal dual N$^*\mathcal{C}$. Then  N$^*\mathcal{C}$   is isomorphic to the normal category $\mathbb{R}(T (\mathcal{C}))$. In particular, the normal dual N$^*\mathcal{C}$ is also a normal category.
\end{thm}

Thus given a regular semigroup $S$, we can associate four normal categories with it: $\mathbb{L}(S)$ and $\mathbb{R}(S)$; their duals N$^*\mathbb{L}(S)$ and N$^*\mathbb{R}(S)$. Their inter relationship can be described using two functors $\Gamma_S\colon  \mathbb{R}(S) \to \text{N}^*\mathbb{L}(S)$ and $\Delta_S\colon  \mathbb{L}(S) \to \text{N}^*\mathbb{R}(S)$ defined as follows: 
\begin{equation}
\begin{split} \label{eqngs}
v\Gamma_S(eS) = H(\rho^e;-) ;\quad& \Gamma_S(\lambda(e,u,f)) = \eta_{\rho^e}\mathbb{L}(S)(\rho(f,u,e),-)\eta_{\rho^f}^{-1};\\
v\Delta_S(Se) = H(\lambda^e;-) \quad& \text{ and }\quad\Delta_S(\rho(e,u,f)) = \eta_{\lambda^e}\mathbb{R}(S)(\lambda(f,u,e),-)\eta_{\lambda^f}^{-1}.
\end{split}
\end{equation}
where $\eta_{\rho^e}$ is the natural isomorphism associated with the $H$-functor $H(\rho^e;-)$ and for $u\in fSe$, the expression $\mathbb{L}(S)(\rho(f,u,e),-)$ is the natural transformation between the covariant hom-functors $\mathbb{L}(S)(Se,-)$ and $\mathbb{L}(S)(Sf,-)$.
To take us to the formal definition of a cross-connection, we shall also require the following notions which abstract the properties of the above defined functors.
\begin{dfn}\label{dfnlociso}
A functor $F$ between two normal categories $\mathcal{C}$ and $\mathcal{D}$ is said to be a \emph{local isomorphism} if $F$ is inclusion preserving, fully faithful and for each $c \in v\mathcal{C}$, $F_{| (  c  ) }$ is an isomorphism of the ideal $ (  c  ) $ onto $  (  F(c)  )  $.
\end{dfn}
Given a normal cone $\gamma \in T (\mathcal{C})$,  the $M$-set of an $H$-functor $H(\gamma;-)$ is defined as:
\begin{equation*}\label{eqnms}
MH(\gamma;-)= \{ c\in v\mathcal{C}: \gamma(c) \text{ is an isomorphism} \}.
\end{equation*}

Observe that for the functors $\Gamma_S$ and $\Delta_S$ defined as above and for objects $Se\in v\mathbb{L}(S)$ and $eS\in v\mathbb{R}(S)$,
\begin{equation*}
Se \in M\Gamma_S(eS) \text{ if and only if } eS\in M\Delta_S(Se).
\end{equation*}
Thus we define:
\begin{dfn} \label{ccxn}
Let $\mathcal{C}$ and $\mathcal{D}$ be normal categories. A \emph{cross-connection} between $\mathcal{C}$ and $\mathcal{D}$ is a quadruple $\Omega=(\mathcal{C},\mathcal{D};{\Gamma},\Delta)$ where $\Gamma\colon  \mathcal{D} \to \text{N}^*\mathcal{C}$ and $\Delta\colon  \mathcal{C} \to \text{N}^*\mathcal{D}$ are local isomorphisms such that for $c \in v\mathcal{C}$ and $d \in v\mathcal{D}$
	\begin{equation*}\label{eqncxnms}
	c \in M\Gamma(d) \iff d\in M\Delta(c).
	\end{equation*}
\end{dfn}

\begin{rmk}	We define a cross-connection using two functors emulating Grillet's \cite{Gri1} original definition using two maps, unlike in \cite{cross0,KSS3} where a cross-connection is defined using a single functor. It can be easily shown that these definitions are equivalent.
\end{rmk}

Summarising the above discussion, Nambooripad described the ideal structure of a regular semigroup using the following theorem.
\begin{thm}\cite[Theorem IV.2]{KSS3}\label{thmcxns}
Let $S$ be a regular semigroup with normal categories $\mathbb{L}(S)$ and $\mathbb{R}(S)$. Define functors $\Gamma_S$ and $\Delta_S$ as in (\ref{eqngs}). Then $\Omega S= (\mathbb{L}(S),\mathbb{R}(S);\Gamma_S,\Delta_S)$ is a cross-connection between  $\mathbb{L}(S)$ and $\mathbb{R}(S)$.
\end{thm}

\medskip

{\bf Cross-connection semigroups}

\medskip

We have seen above how a given regular semigroup induces a cross-connection between its normal categories. Conversely given an abstractly defined cross-connection $\Omega=(\mathcal{C},\mathcal{D};{\Gamma},\Delta)$, Nambooripad gave a recipe to construct its cross-connection semigroup as follows.

Recall that we already have two `generic' regular semigroups: namely $T(\mathcal{C})$ and $T(\mathcal{D})$, the semigroups of normal cones. Then, by the category isomorphism
$$[\mathcal{C},[\mathcal{D},\mathbf{Set}]]\cong [\mathcal{C}\times\mathcal{D},\mathbf{Set}],$$
the functors $\Gamma$ and $\Delta$ induce two bifunctors $\Gamma(-,-)$ and $\Delta(-,-)$, respectively both from $\mathcal{C}\times\mathcal{D}$ to $\mathbf{Set}$. Further, there is a natural isomorphism $\chi: \Gamma(-,-)\to\Delta(-,-)$ between the bifunctors which `links' certain normal cones in $T (\mathcal{C})$ with those in $T (\mathcal{D})$. Finally, the set
\begin{equation}
\mathbb{S}\Omega=\{ (\gamma,\delta) \in T (\mathcal{C})\times T (\mathcal{D}) : (\gamma,\delta) \text{ is linked }\:\}
\end{equation}
is a regular semigroup such that $\mathbb{L}(\mathbb{S}\Omega)$ is isomorphic to $\mathcal{C}$ and $\mathbb{R}(\mathbb{S}\Omega)$ is isomorphic to $\mathcal{D}$. Then, the set
\begin{equation}\label{eqneo}
E_\Omega=\{ (c,d)\in v\mathcal{C}\times v\mathcal{D} : c\in M\Gamma(d)\}
\end{equation}
is the regular biordered set associated with the cross-connection $\Omega$. Here, the pair $(c,d)$ can be identified with a uniquely defined pair of idempotent cones in $T (\mathcal{C})\times T (\mathcal{D})$. Then it can shown that the set $E_\Omega$ is biorder isomorphic to the set $E(\mathbb{S}\Omega)$ of idempotents in $\mathbb{S}\Omega$. The major results \cite[Theorems III.25, IV.17, IV.32 and V.18]{KSS3}  may be summarised as follows.

\begin{thm}\label{summary}
Any regular semigroup $S$ induces a cross-connection between the normal categories $\mathbb{L}(S)$ and $\mathbb{R}(S)$. Conversely any cross-connection $\Omega=(\mathcal{C},\mathcal{D};{\Gamma},\Delta)$ uniquely determines a regular semigroup $\mathbb{S}\Omega$. The category of regular semigroups is equivalent to the category of cross-connections.	
\end{thm}


\section{Outgrowths of Nambooripad's cross-connection theory}

In this section, we briefly record some outgrowths and still ongoing developments of Nambooripad's cross-connection theory.

\medskip

{\bf Inductive groupoids and cross-connections}

\medskip

As discussed in Sections \ref{secind} and \ref{seccxn}, both inductive groupoids and cross-connections are prototypes of regular semigroups. By transitivity, Nambooripad's results (Theorems \ref{thmind} and \ref{summary})  imply that the category of inductive groupoids is equivalent to the category of cross-connections: but these constructions look evidently disconnected. Their interrelationship was explored by Azeef and Volkov \cite{indcxn1,indcxn2} and a direct category equivalence was constructed. In addition to giving the relationship between the ideal structure and the idempotent structure of regular semigroups, their results also give a road map for the transfer of problems of  inductive groupoid origin to the cross-connection framework and vice versa. We proceed to briefly outline this equivalence.

First given a cross-connection $\Omega=(\mathcal{C},\mathcal{D};{\Gamma},\Delta)$, recall that $E_\Omega$ (as defined in equation (\ref{eqneo})) is a regular biordered set. So, we define a category $\mathcal{G}_\Omega$ such that $v\mathcal{G}_\Omega=E_\Omega$ and a morphism in $\mathcal{G}_\Omega$ is a pair of isomorphisms which `respects' the cross-connection. Given a morphism $(f,g)$ from $(c,d)$ to $(c',d')$ and a morphism $(f_1,g_1)$ from $(c_1,d_1)$ to $(c'_1,d'_1)$, define a partial order $\leq_\Omega$ on $\mathcal{G}_\Omega$ as follows:
\begin{multline*}\label{po}
(f,g) \leq_\Omega (f_1,g_1) \iff (c,d) \subseteq (c_1,d_1),\\
(c',d') \subseteq (c'_1,d'_1) \text{ and } (f,g)= ((j(c,c_1)f_1)^\circ,(j(d,d_1)g_1)^\circ).
\end{multline*}
Then it can be easily shown that $ (\mathcal{G}_\Omega,\leq_\Omega)$ forms an ordered groupoid. Further for a suitably defined evaluation functor $\epsilon_\Omega \colon\mathcal{G}(E_{\Omega}) \to \mathcal{G}_\Omega$, it can be verified that the pair $ (\mathcal{G}_\Omega,\epsilon_\Omega)$ forms an inductive groupoid.

Conversely, given an inductive groupoid $(\mathcal{G},\epsilon)$, a cross-connection $(\mathcal{R}_G,\mathcal{L}_G;{\Gamma_G},{\Delta_G})$ was constructed. To build the `left' normal category $\mathcal{L}_G$, three separate categories: one preorder $\mathcal{P}_L$ `responsible' for inclusions, an ordered groupoid $\mathcal{G}_L$ `responsible' for isomorphisms and the last one $\mathcal{Q}_L$ `responsible' for retractions were built and then these categories were `combined' to form $\mathcal{L}_G$. Similarly the normal category $\mathcal{R}_G$ was constructed from three constituent categories and the cross-connection was defined between  $\mathcal{L}_G$ and  $\mathcal{R}_G$. A direct category equivalence between the category of inductive groupoids and cross-connections was also proved \cite[Theorem 5.1]{indcxn2}.

\medskip

{\bf Special classes of Nambooripad's cross-connections}

\medskip

In the first version of cross-connections \cite{cross0}, Nambooripad  discussed several important applications of his theory. We proceed to briefly describe some of them: they give insight as to why Nambooripad pursued such a sophisticated construction.

The first application which Nambooripad considered was the specialisation of his cross-connections to fundamental regular semigroups, thereby extracting Grillet's original theorem regarding {\em fundamental cross-connections} \cite{cross0}. An object $c$ in a normal category $\mathcal{C}$ is said to be {\em reduced} if for morphisms $f,g \in \mathcal{C}(c,c)$, for every $c'\le c$, we have im $j(c',c)f$= im $j(c',c)g$ implies $f=g$.

\begin{dfn}
A normal category $\mathcal{C}$ is said to be  reduced if every object of $\mathcal{C}$ is reduced.
\end{dfn}

Given a regular poset $P$ and $x\in P$, let $P(x) =  \{ y : y\le x \}$. Then we can define a small category $\mathcal{P}$ as follows:
$$v\mathcal{P}=\{P(x): x\in P\} \text{ and } \mathcal{P}(P(x),P(y))=\{f\colon P(x)\to P(y): f\text{ is a normal mapping} \}$$

Then it was shown that the category $\mathcal{P}$ as defined above is a reduced normal category. Also, a  regular semigroup $S$ is fundamental if and only if the normal categories $\mathbb{L}(S)$ and $\mathbb{R}(S)$ are reduced. 
\begin{thm}\cite[Theorem 6.7]{cross0}
Let $S$ be a regular semigroup and let $\Gamma_S$ and $\Delta_S$ be functors as defined in (\ref{eqngs}). Then the pair $(v\Gamma_S,v\Delta_S)$ constitutes a fundamental cross-connection between the regular posets $v\mathbb{L}(S)$ and $v\mathbb{R}(S)$.  Every fundamental cross-connection arises in this way.
\end{thm}


\medskip
In \cite{cross0}, Nambooripad considered the class of strongly regular Baer semigroups and characterised regular semigroups $S$ for which the categories $\mathbb{L}(S)$ and $\mathbb{R}(S)$ have kernels and cokernels.

\begin{thm}
	For a regular monoid $S$ with zero, the following statements are equivalent:
\begin{enumerate}
	\item $S$ is strongly regular Baer semigroup.
	\item $S/\mathscr{L}$ and $S/\mathscr{R}$ are dually isomorphic complemented modular lattices.
	\item The categories $\mathbb{L}(S)$ and  $\mathbb{R}(S)$ have kernels.
	\item The category $\mathbb{L}(S)$ has both kernels and cokernels.
\end{enumerate}
\end{thm}

Further, he introduced the notion of a bounded, pre-additive normal category to study regular rings (in the sense of \cite{good}).
\begin{thm}\cite[Theorem 7.6]{cross0}
A bounded normal category $\mathcal{C}$ is pre-additive if and only if $\mathcal{C}$ is isomorphic to $\mathbb{L}(R)$ for some regular ring $R$.
\end{thm}

Later, Sunny Lukose and Rajan \cite{sunny} studied regular rings 
and obtained an alternate characterisation using $RR$-categories. They showed that the set of normal cones in such a category forms a regular ring. They also described the regular ring of finite rank operators on an infinite dimensional vector space in terms of $RR$-categories.

In \cite{cross0}, Nambooripad also considered {\em semisimple} objects of various types. In particular, he considered a vector space over a field as a $G$-module, for a group $G$ and studied the normal categories and cross-connections which arose in this setting.


\medskip

Several other special classes of the cross-connection construction have been considered in the literature. These works not only clarify the nuances of the complicated construction but also shed light on how the construction could be further generalised.

One of the first studies based on Nambooripad's cross-connection theory was done by Rajendran \cite{rajendthesis,RajNam1,RajNam2} connecting it with bilinear forms. 
A related work was done by Azeef \cite{tlx} from a different perspective wherein the cross-connection structure of several linear transformation semigroups: full, singular and variants, were studied. 
In \cite{tx}, Azeef and Rajan gave the construction of the singular transformation semigroup via cross-connections from the category of subsets and the category of partitions. The variant case was studied in detail in \cite{var} and the cross-connection construction of the regular part of the variant semigroup was provided. 

In \cite{css}, the cross-connection structure of completely simple semigroups were studied and it was shown that the cross-connections are determined by the structure matrices. This result was extended to completely $0$-simple semigroups in \cite{locinverse} by Azeef et al. In addition, the more general class of locally inverse semigroups were discussed in detail in \cite{locinverse}. {\em Unambiguous categories} were introduced as normal categories with unique normal factorisation and unique splitting. It was shown that the category of locally inverse semigroups is equivalent to the category of cross-connections of unambiguous categories.

Several properties of the normal categories of inverse semigroup were identified by Rajan in \cite{invrajan}. A complete characterisation of the inverse case was provided in \cite{locinverse} by Azeef et al. using {\em inversive categories}. Here the structure theorem was obtained using a single category analogous to the ESN Theorem and it was shown that the category of inverse semigroups is equivalent to the category of inversive categories. 

\medskip

{\bf Generalisations of cross-connections}

\medskip

As mentioned earlier, the real inspiration behind the theory of cross-connections is its promise to provide more general structure theorems and its applicability in more general settings. As discussed in Section \ref{secindog}, the inductive groupoid theory was extended in various directions. It must be observed that all these efforts were severely constrained due to the reliance on the set of idempotents of the semigroup.

Nevertheless, as a first step towards extending cross-connection theory to a general setting, it was natural to describe the cross-connection structure of the non-regular classes already studied using the inductive groupoid approach. In this direction, Romeo \cite{romeothesis} studied the cross-connections of concordant semigroups in his PhD thesis under Nambooripad. This work was later refined and published as \cite{conc}. In \cite{conc}, the categories arising from the generalised Green relations of a concordant semigroup were characterised as {\em consistent categories} and a cross-connection was constructed. Conversely given a cross-connection between two consistent categories, a concordant semigroup was obtained as the cross-connection semigroup. This equivalence was shown to be a category equivalence \cite[Theorem 6.8]{conc}. Further, the interrelationship between inductive cancellative categories \cite{armstrong} and cross-connections of consistent categories was studied and a category equivalence was also outlined.

In the late nineties, Nambooripad focussed his attention to a very ambitious problem: to provide a structure theorem for arbitrary semigroups using small categories. He made good progress in this direction \cite{newcross,newcross1}, but unfortunately could not complete the entire construction. We proceed to outline some of his ideas in this problem.

Similar to the regular semigroup case, given an arbitrary semigroup $S$, a category $\mathbb{L}(S)$ is defined as follows:
$$v\mathbb{L}(S) =\{ S^1a: a\in S \}$$
A morphism $\rho \colon S^1a\to S^1b$ is a right translation $x\mapsto xs$ for $s\in S^1$ such that $as=tb$ for some $t\in S^1$. We denote this morphism as $\rho(a,s,b)$. Then it can be seen that given any two morphisms, say $\rho(a,s,b)$ and $\rho(a',s',b')$, they are equal if and only if $a\mathrel{\mathscr{L}}a'$, $b\mathrel{\mathscr{L}}b'$ and $as=as'$.

Observe that a morphism $\rho(a,s,b)$ is an inclusion $j(S^1a,S^1b)$ if and only if $\rho(a,s,b)=\rho(a,1,b)$. These inclusions make  $\mathbb{L}(S)$ a category with subobjects. Then any morphism $\rho(a,s,b)$ can be written as $$\rho(a,s,b)=\rho(a,s,as)\rho(as,1,b)$$ where $\rho(a,s,as)$ is an epimorphism and $\rho(as,1,b)$ is an inclusion. Such a {\em unique} factorisation is called the {\em image factorisation}.

Nambooripad used the natural forgetful functor $U\colon \mathbb{L}(S)\to \mathbf{Set}$ to abstract the above described `set-based' properties of $\mathbb{L}(S)$. This leads to the definition of a {\em set-based category (SBC)} as follows.
\begin{dfn}\label{dfnsbc}
Let	$\mathcal{C}$ be a category with subobjects such that every morphism in $\mathcal{C}$ has a image factorisation. Given a functor $U\colon\mathcal{C}\to \mathbf{Set}$, we say that $\mathcal{C}$ is an {\em SBC} with respect to $U$ if the pair   $(\mathcal{C},U)$ satisfies the following:
\begin{enumerate}
 	\item [(SBC 1)] $U$ is an embedding and $U$ preserves image factorisations.
 	\item [(SBC 1)] For $c,c'\in v\mathcal{C}$ and $x\in U(c)\cap U(d)$, there is a $d\in v\mathcal{C}$ such that
 	$$d\subseteq c,\quad d\subseteq c' \: \text{ and } \: x\in U(d).$$
 \end{enumerate}
\end{dfn}

In this setting, Nambooripad proved some preliminary results and further conjectured that  an appropriately defined cross-connection between two SBCs will give rise to a semigroup. Conversely, any arbitrary semigroup $S$ determines a cross-connection between its constituent SBCs $\mathbb{L}(S)$ and $\mathbb{R}(S)$ such that its cross-connection semigroup provides a natural representation of $S$. This remains a major open problem.



\section{Connections between regular semigroups and other areas}

In addition to his deep work on the structure of regular semigroups, Nambooripad had  interests in many areas of mathematics and in particular in exploring connections between regular semigroups and other areas. He wrote several papers along these lines, some in collaboration with some of his PhD students, and several of his PhD students worked on such connections.  In particular, his work revealed interesting connections between regular semigroups and linear algebra (stochastic matrices  \cite{premthesis}, singular matrices \cite{geethathesis,GeeNam},  the geometry and topology of idempotent matrices
\cite{vnkthesis,KrishnaNam} and bilinear forms \cite{rajendthesis,RajNam1,RajNam2}). He also studied group actions on lattices \cite{KSS9,KSS10}. He had a particular interest in connections between regular semigroups and operator algebras. We briefly summarize some of his work and that of some of his students in operator algebras, in particular in the study of Fredholm operators, finite rank operators on Hilbert space, and von Neumann algebras.
 Much of his work in this direction is incomplete and suggests that it may be very fruitful to pursue additional work in this direction.

\medskip

{\bf The semigroup of Fredholm operators }

\medskip

In his PhD thesis \cite{ekthesis}, E. Krishnan studied the category of Fredholm operators between
topological vector spaces and the  semigroup of Fredholm operators on such a  space. Much
 of his work is contained in his joint paper with Nambooripad \cite{KrishNam}. We briefly describe some of the ideas discussed in this paper and refer the reader to the original paper for more detail. All spaces under consideration will be Hausdorff locally convex topological vector spaces with underlying field either the real or complex numbers.

A {\em Fredholm operator} $f : X \to Y$ between topological vector spaces is a continuous linear map that is an open map onto its range $im(f)$ such that $im(f)$ is a closed subspace of finite codimension in $Y$ and the null space $N(f)$ is a finite dimensional subspace of $X$.

Among many other things, Krishnan and Nambooripad \cite[Proposition 3.21]{KrishNam} prove that if $f : X \to Y$ and $g: Y \to Z$ are Fredholm operators then their composition $gf : X \to Z$ is also a Fredholm operator, so the class  of locally convex spaces together with the Fredholm operators between
them forms a category, denoted by $\mathcal{F}$.

In order to study the category $\mathcal{F}$ in more detail, Krishnan and Nambooripad \cite{KrishNam} introduce a theory of {\em regular} categories, along the lines of Nambooripad's theory of regular semigroups. Here a morphism $f \in \mathcal{C}(v,w) $ in a  category $\mathcal C$ is called {\em regular} if there is a morphism $f' \in \mathcal{C}(w,v)$ in $\mathcal C$ with $f = ff'f$. We refer the reader to Krishnan and Nambooripad's paper \cite{KrishNam} for full details of the definition and basic properties of regular categories. Among other things, they show that  every Fredholm operator $f$ in $\mathcal F$ is a {\em regular morphism} in this category and hence that, for each locally convex space $X$,  the set $\mathcal F (X)$ of Fredholm operators with domain and codomain equal to $X$ forms a regular semigroup \cite[Proposition 3.23]{KrishNam}.

Krishnan and Nambooripad \cite{KrishNam} prove that, with respect to the topology of uniform convergence on bounded sets, the semigroup $\mathcal{F}(X)$ is a semitopological semigroup (that is, the multiplication is separately continuous) for any locally convex space $X$, and is in fact a topological semigroup (that is, the multiplication is jointly continuous) if $X$ is a normed space. The paper \cite{KrishNam} studies both the  topological and algebraic properties of the semigroup $\mathcal{F}(X)$ in some detail. In particular, they prove \cite[Theorem 4.1]{KrishNam} that if $X$ is a locally convex space, then the relation

\[\mu = \{(f,g) \in \mathcal{F} \times \mathcal{F} : g = \lambda f \text{ for some scalar }\lambda \neq 0\}\]

\noindent is the maximum idempotent-separating congruence on $\mathcal{F}(X)$. They are able to use this to study the fundamental representation of $\mathcal{F}(X)$: in particular, by making use of some of the results of Nambooripad and Pastijn \cite{NamPas2}, they show that if $X$ is an infinite-dimensional locally convex space, then the elements of $\mathcal{F}(X)/{\mu}$ may be identified with projective maps on the projective geometry of subspaces of $X$. Much additional information about several important congruences on $\mathcal{F}(X)$ is provided in \cite{KrishNam}.

Krishnan and Nambooripad also introduce a new integral invariant $k(X)$ of a topological vector space $X$ as follows. For each locally convex Hausdorff space $X$ they define $k(X)$ to be the minimum codimension of a proper closed subspace $Y$ of finite codimension in $X$ such that $Y \cong X$, if such a subspace $Y$ exists, and $k(X) = 0$ if no such subspace exists. They make use of this invariant to study several algebraic properties of the semigroup $\mathcal{F}(X)$. For example, they prove the following theorem \cite[Theorems 5.9, 5.11 and 5.12]{KrishNam}.

\begin{thm}
\label{KNthm}
Let $X$ be a locally convex Hausdorff space. Then

(a) $\mathcal{F}(X)$ is completely semisimple and unit regular if and only if $k(X) = 0$.

(b) $\mathcal{F}(X)$ is bisimple if and only if $k(X) = 1$.

(c) $\mathcal{F}(X)$ is simple if and only if $k(X) \neq 0$. (So $\mathcal{F}(X)$ is simple but not bisimple if and only if $k(X) > 1$.)

\end{thm}

The paper \cite{KrishNam} and Krishnan's thesis \cite{ekthesis} contain a wealth of additional information about Fredholm operators from the point of view of Nambooripad's theory of regular semigroups.

\medskip

{\bf Finite rank operators on Hilbert space}

\medskip

In her thesis \cite{shlthesis} completed under the direction of K.S.S. Nambooripad, Sherly Valanthara studied the semigroup of finite rank bounded operators on a  Hilbert space. Her thesis contains a wealth of information about the algebraic and topological properties of the semigroup $\mathcal{B}(\bf H)$ of continuous (i.e. bounded) operators on a Hilbert space $\bf H$ and certain of its subsemigroups. Some of her results are contained in her joint papers with Krishnan \cite{eksh1,eksh}. She provides a very nice self-contained account of several established and new algebraic and topological properties of regular elements, idempotents and Moore-Penrose inverses of elements in $\mathcal{B}(\bf H)$ and a detailed study of the structure of the semigroup $\mathcal{K}(\bf H)$ of finite rank operators in $\mathcal{B}(\bf H)$.  We collect several of her results \cite[Propositions 3.1.2, 3.1.3, 3.2.6, 3.3.2, 3.3.3, 3.3.7, 3.3.8]{shlthesis} about the algebraic properties of $\mathcal{K}(\bf H)$ in the following theorem.

\begin{thm}
\label{finiterankalg}

Let $\mathcal{K}(\bf H)$ denote the semigroup of finite rank continuous operators on a Hilbert space $\bf H$. Then

(a) $\mathcal{K}(\bf H)$ is a self-adjoint, regular, completely semisimple subsemigroup of $\mathcal{B}(\bf H)$.

(b) $\mathcal{K}(\bf H)$ is idempotent-generated if $\bf H$ is infinite-dimensional.

(c) If $f$ is an operator in $\mathcal{K}(\bf H)$ then its Moore-Penrose inverse $f^{\dagger}$ is also in $\mathcal{K}(\bf H)$.

(d) If $e$ and $f$ are idempotents in $\mathcal{K}(\bf H)$ with $e \mathcal{D} f$, then $e$ and $f$ are connected by an $E$-chain in $E(\mathcal{K}(\bf H))$ of length at most $3$.

(e) The set $\mathcal{K}(\bf H) \cup \mathcal{G}(\bf H)$ is a strongly unit regular subsemigroup of $\mathcal{B}(\bf H)$. (Here $\mathcal{G}(\bf H)$ is the group of invertible operators in $\mathcal{B}(\bf H)$.)

\end{thm}
We remark that part (d) of Theorem \ref{finiterankalg} is an extension to the setting of finite-rank operators of infinite-dimensional Hilbert spaces of the corresponding result for the multiplicative semigroup $M_n(F)$ of $n \times n$ matrices over a field $F$, due independently to Pastijn \cite{Pas2} and Putcha \cite{Putcha2} as a result of their constructions of the biordered set of $M_n(F)$.

With respect to the topology induced by the operator norm, the semigroup $\mathcal{B}(\bf H)$ and all its subsemigroups are  topological semigroups. Sherly Valanthara provides much information about the topological properties of this semigroup in her thesis. In particular, the following theorem collects some of these properties from her thesis \cite[Propositions 3.4.1, 3.4.2, 3.4.7, 3.4.9, 3.4.10]{shlthesis}.

\begin{thm}
\label{finiteranktop}

(a) The semigroup $\mathcal{K}(\bf H)$ is a topological semigroup with continuous involution $f \mapsto f^*$. Furthermore, the map $f \mapsto f^{\dagger}$ is continuous on each non-zero $\mathcal H$-class of $\mathcal{K}(\bf H)$.

(b) Each  element of $\mathcal{K}(\bf H)$ is the limit of a sequence of invertible operators on $\mathcal{B}(\bf H)$.

(c) A pair of idempotents in $\mathcal{K}(\bf H)$ are connected by a path in the topological space $E(\mathcal{K}(\bf H))$ if and only if they are connected by an $E$-sequence in the biordered set $E(\mathcal{K}(\bf H))$.

(d) The path components of the space $E(\mathcal{K}(\bf H))$ are the sets of idempotents in the various $\mathcal{D}$-classes of $\mathcal{K}(\bf H)$.

\end{thm}

The thesis \cite{shlthesis} is full of many additional interesting results connecting the algebraic and topological properties of the semigroup $\mathcal{K}(\bf H)$ with Nambooripad's work on the structure of regular semigroups.

\medskip

{\bf von Neumann algebras}

\medskip

Nambooripad was very interested in making use of the theory of regular semigroups to obtain some information about von Neumann algebras. His paper \cite{KSS7} and a personal communication \cite{KSS8} provides some information along these lines. We provide a brief discussion of some of his ideas about connecting regular semigroups and von Neumann algebras.

The theory of von Neumann algebras traces its origins back to the original work of von Neumann in 1930 and a series of papers by Murray and von Neumann in the 1930's. A {\em von Neumann algebra} is a $^*$-algebra of bounded operators on a Hilbert space that contains the identity operator and is closed in the weak operator topology. Equivalently, a von Neumann algebra may be defined as a subset of the set $\mathcal{B}({\bf H})$ of bounded operators on a Hilbert space ${\bf H}$ that is closed under the $^*$-operation and is equal to its double commutant. (The equivalence of these definitions is von Neumann's ``double commutant" theorem). We refer the reader to the books by Fillmore \cite{Fill} or Sunder \cite{Sund} or Dixmier \cite{Dix} for this result and the basic notation, concepts, references,  and many other standard results about von Neumann algebras.

Let $\mathcal{A} \subseteq \mathcal{B}({\bf H})$ be a von Neumann algebra acting on the Hilbert space $\bf H$ and let $P(\mathcal{A}) = \{p \in \mathcal{A} : p = p^* = p^2\}$ be the set of {\em projections} of $\mathcal{A}$. The range $im(p)$ of a projection is a closed subspace of $\bf H$. The projections in $\mathcal{A}$ are  the operators that give an orthogonal projection onto some closed subspace of  $\bf H$.  It is well-known \cite{Dix,Sund,Fill} that $\mathcal{A}$ is generated in $\mathcal{B}({\bf H})$ by $P(\mathcal{A})$: that is, a von Neumann algebra is uniquely determined by its projections and the underlying Hilbert space on which it acts. (Here by the von Neumann algebra {\em generated} by a subset $A \subseteq \mathcal{B}({\bf H})$ we mean the double commutant $(A \cup A^*)''$, that is, the weak closure of the $^*$-algebra generated  by $1$ and $A$ in $\mathcal{B}({\bf H})$).  Since projections are in particular idempotents, it follows that $\mathcal{A}$ is uniquely determined in the sense described above by the biordered set $E(\mathcal{A})$ of idempotents of the multiplicative semigroup of $\mathcal{A}$.

These observations led Nambooripad to study the structure of the biordered set $E(\mathcal{A})$ of a von Neumann algebra $\mathcal{A}$. Recall that if $p$ is a projection in $\mathcal{A}$ then its range $im(p)$ is a closed linear subspace of $\bf H$. A subspace of $\bf H$ is said to {\em belong to $\mathcal{A}$}  if it is the image of some projection in $\mathcal{A}$. The map $\theta: p \mapsto im(p)$ is a one-one correspondence between the set of projections of $\mathcal{A}$ and the set $\Theta(\mathcal{A})$ of closed subspaces of $\bf H$ that belong to $\mathcal{A}$. In fact the partially ordered set of projections of a von Neumann algebra $\mathcal{A}$ and the partially ordered set $\Theta(\mathcal{A})$ of subspaces of $\bf H$ belonging to $\mathcal{A}$ (ordered with respect to set inclusion) are isomorphic {\em continuous geometries} in the sense of von Neumann \cite{vonN}. In particular, $P(\mathcal{A})$ and $\Theta(\mathcal{A})$ are {\em complemented modular lattices}. Nambooripad made use of this, together with Pastijn's construction of biordered sets from complemented modular lattices \cite{Pas2} to study the biordered set $E(\mathcal{A})$. Let

\[E'(\mathcal{A}) = \{(N,V): N,V \in \Theta(\mathcal{A}), \, N \wedge V = \{0\} \text{ and }N \vee V = {\bf H}\}.\]

 We include an outline of Nambooripad's proof  \cite{KSS8} of the following fact since he did not publish the result as far as we are aware.

\begin{thm}
\label{vnbiorder}
Let $\mathcal{A}$ be a von Neumann algebra acting on a Hilbert space $\bf H$. Then the set $E'(\mathcal{A})$ admits the structure of a regular biordered set and the biordered set $E(\mathcal{A})$ is isomorphic to a biordered subset of $E'(\mathcal{A})$.

\end{thm}

\noindent {\bf Outline of Nambooripad's proof. }

By Theorem 1 of Pastijn's paper \cite{Pas2}, $E'(\mathcal{A})$ is a biordered set with quasi-orders and basic products determined as follows:

\[(M,U) \omega^l(N,V) \text{ iff } U \subseteq V \text{ and in this case } (N,V)(M,U) = (N \vee (V \wedge M),U)\]
and
\[(M,U) \omega^r(N,V) \text{ iff } N \subseteq M \text{ and in this case } (M,U)(N,V) = (M,V \wedge (N \vee U))\]

For $(L,U),(M,V) \in E'(\mathcal{A})$, the sandwich set $\mathcal{S}((L,U),(M,V))$ consists of all elements $(N,W) \in E'(\mathcal{A})$ such that $N$ is any complement of $U \vee M$ in the interval $[M,1]$ and $W$ is any complement of $U \wedge M$ in the interval $[0,U]$. Since these relative complements exist, all sandwich sets are non-empty and so $E'(\mathcal{A})$ is a regular biordered set.

Now let $e$ be any idempotent in $E(\mathcal{A})$. Then the range $im(e)$ and null space ${\bf N}(e)$ are complementary subspaces of $\bf H$. The fact that $e$ is a regular element of $\mathcal{B}(\bf H)$, implies that $im(e)$  is a closed subspace of  $\bf H$ and that $ee^{\dagger}$ is a projection onto $im(e)$, where $e^{\dagger}$ is the Moore-Penrose inverse of $e$. By Proposition 9 of Nambooripad's paper \cite{KSS7}, $e^{\dagger}$ is in $\mathcal{A}$.   Hence  $im(e) = im(ee^{\dagger})$ belongs to $\mathcal{A}$. A routine calculation shows that $im(1 - e) = {\bf N}(e)$ and $(1-e)$ is an idempotent in $\mathcal{A}$ so   ${\bf N}(e)$ is also a closed subspace of $\bf H$ belonging to $\mathcal{A}$. So the pair $({\bf N}(e),im(e))$ is in $E'(\mathcal{A})$.
Hence the map $\epsilon: e \mapsto ({\bf N}(e),im(e))$ defines an injection from $E(\mathcal{A})$ into
$E'(\mathcal{A})$. It is routine to check that $\epsilon$ preserves basic products, so $E(\mathcal{A})$ is isomorphic to a biordered subset of $E'(\mathcal{A})$.  $\Box$

\medskip

Nambooripad also obtained some interesting information about regular elements in von Neumann algebras. It is well known that an element $\gamma$ in $\mathcal{B}({\bf H})$ is regular if and only if its range $im(\gamma)$ is closed, and in this case ${\gamma}{\gamma}^{\dagger}$ is the projection on $im(\gamma)$ and ${\gamma}^{\dagger}\gamma$ is the projection on ${\bf N}(\gamma)^{\perp}$. The multiplicative semigroup of $\mathcal{B}({\bf H})$ is not regular if $\bf H$ is infinite-dimensional and so the multiplicative subsemigroup of a von Neumann algebra $\mathcal{A}$ is not in general regular, and in fact the subsemigroup of $\mathcal{A}$ generated by the regular elements is not a regular semigroup. However, for factors of type $I_n$ or $II_1$ (see \cite{Sund} for precise definitions of these classes of von Neumann algebras), the projection lattice is modular. This enables Nambooripad to prove the following theorem \cite[Proposition 11]{KSS7}.

\begin{thm}
\label{typereg}
If $\mathcal{A}$ is a von Neumann algebra that is a factor of type $I_n$ or $II_1$ then the set of regular operators in $\mathcal{A}$ forms a regular subsemigroup of $\mathcal{A}$.

\end{thm}
From Theorem \ref{typereg} it is clear that if $\mathcal{A}$ is a von Neumann algebra of type $I_n$ of $II_1$, then $E(\mathcal{A})$ is a regular biorderd set. It is not known whether in fact $E(\mathcal{A})$ is a regular biordered set for some larger class of von Neumann algebras. It appears that additional work exploring connections between Nambooripad's theory of regular semigroups and the structure of von Neumann algebras may prove fruitful.

\medskip


\section{PhD theses directed by Nambooripad}

The topics in which Nambooripad directed PhD theses arise from a
wide range of areas in mathematics. In particular, he directed theses in the study of semigroup-theoretic aspects  of  matrix theory, the theory of operators on  Hilbert spaces, and the geometry and toplogy of linear spaces.

The use of category theory in the presentation of interrelationships between
various classes has been a favourite style in Nambooripad's works.
His  cross connection theory for the structure of  regular semigroups
can be seen as an instance of heavy use of category theoretic ideas.

The PhD thesis of A.R. Rajan\cite{Rajan1} on combinatorial regular semigroups
and the cross connection related theses of D. Rajendran \cite{rajendthesis} and
P.G. Romeo \cite{romeothesis} use categories as a major tool in the
formulation of their results.

S. Premchand's thesis \cite{premthesis} on stochastic matrices, K. Geetha's
 thesis \cite{geethathesis} on singular matrices and V.N. Krishnachandran's thesis \cite{vnkthesis} on the geometry and topology
of idempotent matrices  provide a deep study of matrix theory using semigroup
theoretic tools.

E. Krishnan's thesis \cite{ekthesis} on Fredholm operators and Sherly
Valanthara's thesis \cite{shlthesis}  on
finite rank operators contain a significant amount of operator theory on
Hilbert spaces.

R. Veeramony's  thesis \cite{rvthesis}  on subdirect products
and S. Radhakrishnan Chettiar's thesis \cite{srcthesis} on extensions of regular semigroups
provide a deep study of some theoretical aspects of semigroups.

The following is the list of PhD's directed by Nambooripad.
All of these were awarded from the University of Kerala.
\begin{enumerate}
	\item R. Veeramony: {\em Subdirect products of regular semigroups} (1981)
	\item A.R. Rajan: {\em Structure of combinatorial regular semigroups} (1981)
	\item S. Premchand: {\em Semigroup of stochastic matrices} (1985)
	\item E. Krishnan: {\em The semigroup of Fredholm operators} (1990)
	\item P.G. Romeo: {\em Cross connections of  concordant semigroups} (1993)
	\item D. Rajendran: {\em Cross connection of linear transformation semigroups} (1995)
	\item K. Geetha: {\em Semigroup of singular matrices} (1995)
	\item S. Radhakrishnan Chettiar: {\em A study on extensions of  regular semigroups} (1996)
	\item V. N. Krishnachandran: {\em The topology and geometry of the biordered set of
	idempotent matrices} (2001)
	\item Sherly Valanthara: {\em The semigroup of finite rank operators} (2002)
\end{enumerate}



\end{document}